\def\N{\mathbb N}
\def\Z{\mathbb Z}
\def\R{\mathbb R}
\def\Q{\mathbb Q}

\def\ol{\overline}
\def\pfz{\begin{proof}}
\def\pfk{\end{proof}}

\documentclass{article}
\usepackage{amssymb,amsmath,amsthm}
\usepackage{rotating}
\usepackage{algorithmic,algorithm}
\usepackage{lscape}
\newtheorem{lem}{Lemma}
\newtheorem{thm}[lem]{Theorem}
\newtheorem{prop}[lem]{Proposition}
\newtheorem{coro}[lem]{Corollary}

\newtheorem{ex}[lem]{Example}

\begin{document}
\title{Arithmetics in number systems with a negative base}
\author{Z. Mas\'akov\'a\footnote{e-mail: zuzana.masakova@fjfi.cvut.cz}, E. Pelantov\'a, T. V\'avra\\[1mm]
{\normalsize Department of Mathematics FNSPE, Czech Technical University in Prague}\\
{\normalsize Trojanova 13, 120 00 Praha 2, Czech Republic}}
\date{}

\maketitle

\begin{abstract}
We study the numeration system with a negative base, introduced by
Ito and Sadahiro. We focus on arithmetic operations in the set
${\rm Fin}(-\beta)$ and $\Z_{-\beta}$ of numbers having finite
resp.\ integer $(-\beta)$-expansions. We show that ${\rm
Fin}(-\beta)$ is trivial if $\beta$ is smaller than the golden
ratio $\frac12(1+\sqrt5)$. For $\beta\geq\frac12(1+\sqrt5)$ we
prove that ${\rm Fin}(-\beta)$ is a ring, only if $\beta$ is a
Pisot or Salem number with no negative conjugates. We prove the
conjecture of Ito and Sadahiro that ${\rm Fin}(-\beta)$ is a ring
if $\beta$ is a quadratic Pisot number with positive conjugate.
For quadratic Pisot units we determine the number of fractional
digits that may appear when adding or multiplying two
$(-\beta)$-integers.
\end{abstract}

\section{Introduction}

There exist many ways to represent real numbers. Besides commonly
used decimal and binary number system, there is for example the
well known expansion of numbers in the form of a continued
fraction. The choice of number representation depends on the
purpose for which it will be used. The conventional representation
of real numbers in base $b\in\N$, $b\geq 2$, has been generalized
in 1957 by A. R\'enyi~\cite{Renyi} when he considered writing
non-negative real $x$ in the form
\begin{equation}\label{eq:1}
x=\sum_{i=-\infty}^{k} x_i\beta^i\,,\qquad x_i\in\{0,1,\dots,\lceil\beta\rceil-1\}\,,
\end{equation}
where for the base $\beta$ one chooses arbitrary real $\beta>1$.
The choice of an irrational base brings into play new phenomena
which are not found in the numeration systems with integer base.
For example, the set ${\rm Fin}(\beta)$ of such real $x$ that have
only a finite number of non-zero digits $x_i$ in the
expression~\eqref{eq:1} need not to be closed under addition. An
important motivation for the study of R\'enyi $\beta$-expansions was
given by mathematical models of non-crystalline solids with
long-range order, the so-called quasicrystals. It turned out that
a suitable discrete set for labelling coordinates of atoms in
quasicrystals is formed by the $\beta$-integers. These are real
numbers $x$ which can be written as $\sum_{i=0}^k x_i\beta^i$,
where the choice of base $\beta$ is related to the rotational
symmetry displayed by the material. Many papers since then were
devoted to the study of properties of the set ${\rm Fin}(\beta)$
and $\beta$-integers $\Z_\beta$, to algorithms for addition and
multiplication, their relation to aperiodic tilings of the space,
etc.

Recently, Ito and Sadahiro~\cite{ItoSadahiro} have suggested to
study representations of real numbers in base $-\beta$, where
$\beta>1$, i.e.\ in the form
$$
x=\sum_{i=-\infty}^{k} x_i(-\beta)^i\,,\qquad x_i\in\{0,1,\dots,\lfloor\beta\rfloor\}\,.
$$
The advantage of this numeration system with negative base stems
from the fact that both positive and negative real numbers can be
represented with non-negative digits, without the necessity to
mark the sign of the number. In their paper, Ito and Sadahiro
derive basic properties of the digit strings
$x_kx_{k-1}x_{k-2}\cdots$ corresponding to $(-\beta)$-expansions
of real numbers $x$.

The aim of this article is to study arithmetical properties of the
numeration system with negative base. The paper is organized as
follows: In Section~\ref{sec:prelim} we recall the definition of
R\'enyi $\beta$-expansions, and several number-theoretical notions
which allow us to state known results about sets ${\rm Fin}(\beta)$
and $\Z_\beta$. We give the definition of $(-\beta)$-expansions
and cite relevant facts from the paper~\cite{ItoSadahiro}. In
Section~\ref{sec:-beta} we study the behaviour of the set ${\rm
Fin}(-\beta)$ and $\Z_{-\beta}$ of finite resp.\ integer
$(-\beta)$-expansions with respect to arithmetic operations.
Section~\ref{sec:kvadratic} is devoted to the properties of ${\rm
Fin}(-\beta)$ when $\beta$ is a quadratic number. We prove the
conjecture of Ito and Sadahiro~\cite{ItoSadahiro} that ${\rm
Fin}(-\beta)$ is a ring if $\beta$ is a quadratic Pisot number
with a positive conjugate. For quadratic Pisot units we determine
the number of fractional digits that may appear when adding or
multiplying two $(-\beta)$-integers.

\section{R\'enyi expansions}\label{sec:prelim}

Let us recall representations of numbers in the numeration system
with a positive real base, as introduced by R\'enyi~\cite{Renyi}.
Let $\beta>1$. For a non-negative real number $x$, the well known
greedy algorithm yields a representation of $x$ in the form
$$
x=\sum_{i=-\infty}^k x_i\beta^i\,,
$$
where the digits $x_i$ are integers
$\{0,1,\dots,\lceil\beta\rceil-1\}$. Such representation is called
the $\beta$-expansion of $x$. If $\beta\notin\N$, then not all
digit strings can arise as $\beta$-expansion of some real number
$x$. In order to describe the digit strings that are admissible as
$\beta$-expansions, one defines the so-called R\'enyi expansion of
1, denoted by $d_\beta(1)$. It is a sequence
$d_\beta(1)=t_1t_2t_3\cdots$ of digits in
$\{0,1,\dots,\lceil\beta\rceil-1\}$, such that
$t_1=\lfloor\beta\rfloor$ and $\sum_{i=1}^\infty t_{i+1}\beta^i$
is a $\beta$-expansion of $\beta-\lfloor\beta\rfloor$. A real
number $\beta>1$ is called Parry, if the R\'enyi expansion of 1 is
an eventually periodic sequence, $d_\beta(1)=t_1\cdots
t_m(t_{m+1}\cdots t_{m+p})^\omega$, $t_m\neq t_{m+p}$. Note that $d_\beta(1)$ is never purely periodic.
We distinguish simple Parry numbers as those with period formed only
by 0's. We define
$$
d^*_\beta(1)=\left\{\begin{array}{ll}
\big(t_1\cdots t_{m-1}(t_{m}-1)\big)^\omega &\hbox{ if } d_\beta(1)=t_1\cdots t_m0^\omega\,,\ t_m\neq0\,,\\[2mm]
d_\beta(1)&\hbox{ otherwise.}
\end{array}\right.
$$
Parry~\cite{Parry} has shown that a digit string $x_kx_{k-1}\cdots$ is admissible if and only if each of its suffixes satisfies
$$
x_ix_{i-1}\cdots \prec d^*_\beta(1)\,,
$$
where $\prec$ is the lexicographic order on strings.

The lexicographic order on the digit strings obtained by the
greedy algorithm corresponds to the natural order on the real
line. More precisely, if $x=\sum_{i=-\infty}^kx_i\beta^i$,
$y=\sum_{i=-\infty}^ky_i\beta^i$ are the $\beta$-expansions of $x$
and $y$, respectively, then
$$
x<y \quad\iff\quad x_kx_{k-1}x_{k-2}\cdots \prec y_ky_{k-1}y_{k-2}\cdots\,.
$$

Several classes of numbers display exceptional properties when
taken as bases $\beta$ of R\'enyi numeration systems. Except Parry
numbers, these are namely Pisot numbers and Salem numbers. A real
number $\beta>1$ is called Pisot, if is an algebraic integer with
all conjugate in the interior of the unit disc. $\beta$ is a Salem
number, if all conjugates lie inside the unit disc and at least
one on the unit circle. It can be shown that all Pisot numbers are
Parry, the same statement for Salem numbers is an unproved
conjecture.

When studying arithmetical properties of the $\beta$-numeration
system, one is interested in the set of finite $\beta$-expansions,
denoted by ${\rm Fin}(\beta)$. Question arises, for which bases
$\beta$ this set is closed under addition, subtraction and
multiplication, i.e.\ has a ring structure. Such bases are said to
have the finiteness property. Frougny and Solomyak~\cite{FruSo}
have shown that a necessary condition for $\beta$ to have the
finiteness property is to be a Pisot number. The converse is not
true. Akiyama~\cite{akiyamaCubic} has described all cubic Pisot
units for which ${\rm Fin}(\beta)$ is a ring, other results about
this problem are found in~\cite{FruSo,Hollander}, and its
connection to tilings is explained in~\cite{akiyamaTiles}. Another
question to ask is about the time and space complexity of the
arithmetical operations over finite $\beta$-expansions. One
measure for this are the values $L_\oplus(\beta)$,
$L_\otimes(\beta)$ denoting the maximal length of the
$\beta$-fractional part arising in addition, resp. multiplication
of numbers. More precisely, denoting by $\Z_\beta$ the set of
numbers whose $\beta$-expansion has non-zero digits only for
non-negative indices, $$ \Z_\beta = \Big\{\pm x\,\Big|\,
x=\sum_{i=0}^k x_i\beta^i\ \hbox{is the $\beta$-expansion of
$x$}\Big\}\,,
$$
one puts
\begin{align*}
L_\oplus(\beta)&=\min\{l\in\N\mid \forall\,x,y\in\Z_\beta,\ x+y\in{\rm Fin}(\beta)\Rightarrow x+y\in\beta^{-l}\Z_\beta\}\,,\\
L_\otimes(\beta)&=\min\{l\in\N\mid \forall\,x,y\in\Z_\beta,\ x\cdot y\in{\rm Fin}(\beta)\Rightarrow x\cdot y\in\beta^{-l}\Z_\beta\}\,.
\end{align*}
Bernat~\cite{Bernat} shows that the
values $L_\oplus(\beta)$, $L_\otimes(\beta)$ are finite for every Perron number, i.e.\ for every algebraic integer $\beta>1$
with conjugates in modulus smaller than $\beta$.
Papers~\cite{BernatCubic,BuFrGaKr,GuMaPe,AFMP} provide bounds on
$L_\oplus(\beta)$, $L_\otimes(\beta)$ for specific classes of
Pisot numbers.

%
%
%

\section{Ito-Sadahiro expansions}

In this paper we study arithmetical properties of the $(-\beta)$-numeration system, introduced by Ito and Sadahiro~\cite{ItoSadahiro}.
Let $\beta>1$. Any real number $x$ can be expressed in the form
\begin{equation}\label{eq:vyjadreni}
x=\sum_{i=-\infty}^k x_i(-\beta)^i\,,\qquad x_i\in\Z\,.
\end{equation}
Symbolically, the number $x$ can be written as a sequence of digits $x_i$, where the delimiter ${\scriptsize \bullet}$
separates between coefficients at non-negative and negative powers of $(-\beta)$, i.e.
$$
x_kx_{k-1}\cdots x_1x_0{\scriptsize \bullet}x_{-1}x_{-2}\cdots \hbox{ if $k\geq0$,}\quad\hbox{or}\quad 0{\scriptsize \bullet}0^{-k-1}x_{k}x_{k-1}\cdots
\hbox{ if $k<0$.}
$$
Ito and Sadahiro give a prescription to obtain a $(-\beta)$-representation of numbers $x\in I_\beta=\big[\frac{-\beta}{\beta+1},\frac1{\beta+1}\big)$
using the transformation $T_{-\beta}:I_\beta\mapsto I_\beta$,
$$
T_{-\beta}(x)=-\beta x-\Big\lfloor-\beta x+\frac{\beta}{\beta+1}\Big\rfloor\,,
$$
and set
$$
x_i=\big\lfloor-\beta T^{i-1}_{-\beta}(x)+\frac{\beta}{\beta+1}\big\rfloor\qquad\hbox{and}\qquad d_{-\beta}(x)=x_1x_2x_3\cdots\,.
$$
As shown in~\cite{ItoSadahiro}, a very important role is played by the digit string
$$
d_{-\beta}(\ell_\beta) = d_1d_2d_3\cdots\,, \qquad \hbox{where }\ \ell_{\beta}=\frac{-\beta}{\beta+1}\,.
$$

Obviously,  $0{\scriptsize \bullet} x_1x_2x_3\cdots$ is a $(-\beta)$-representation of $x$ for any $x\in I_\beta$.
The transformation $T_{-\beta}$ can be used to find a $(-\beta)$-representation of any real number $x$, by using $d_{-\beta}\big(x(-\beta)^{-j}\big)$ for a suitable $j\in\N$ and by suitable placement of the delimiter ${\scriptsize \bullet}$.
However, the $(-\beta)$-representation obtained by such a procedure depends on the choice of $j$, since for the left end-point $\ell_{\beta}=\frac{-\beta}{\beta+1}$ of the interval $I_\beta$,
$$
d_{-\beta}\big(\ell_\beta\big)=d_1d_2d_3\cdots\quad\hbox{ and }\quad d_{-\beta}\big((-\beta)^{-2}\ell_\beta\big)=1d_1d_2d_3\cdots\,.
$$
In order to define a unique $(-\beta)$-expansion for every real number $x$, we choose the representation which satisfies
the natural property of usual number representations that multiplication by the base results in shifting the digit sequence. More formally,
we require that for all $x\in\R$ and for all $j\in\Z$,
\begin{equation}\label{eq:ekvivalence}
\begin{aligned}
x&=\sum_{i=-\infty}^k x_i(-\beta)^i\ \hbox{is the $(-\beta)$-expansion of $x$ }  \\
&\hspace*{2.5cm}\quad\Updownarrow\quad\\
(-\beta)^j x&=\sum_{i=-\infty}^{k+j} x_{i-j}(-\beta)^i\ \hbox{is the $(-\beta)$-expansion of $(-\beta)^jx$.}
\end{aligned}
\end{equation}
Among all $(-\beta)$-representations of $x$, we define
the $(-\beta)$-expansion of $x$ as the one which is obtained by the following algorithm.


\smallskip\noindent
\parbox{\textwidth}{
\hrule
  \begin{algorithmic}
    \REQUIRE $x\in\R$
    \STATE Find the minimal non-negative integer $j$ such that $y=x(-\beta)^{-j}\in \big(\frac{-\beta}{\beta+1},\frac1{\beta+1}\big)$.
    \STATE Find $d_{-\beta}(y)=y_1y_2y_3\cdots$
    \IF{$j=0$}
    \STATE put $\langle x\rangle_{-\beta}=0{\scriptsize \bullet}y_1y_2y_3\cdots$.
    \ELSE
    \STATE put $\langle x\rangle_{-\beta}=y_1\cdots y_{j}{\scriptsize \bullet}y_{j+1}y_{j+2}\cdots$.
    \ENDIF
   \end{algorithmic}
\hrule}
\smallskip

Applying the algorithm to the left end-point $\ell_{\beta}$ of the interval $I_\beta$, we obtain
$$
\langle \ell_\beta\rangle_{-\beta} = 1d_1{\scriptsize \bullet} d_2d_3\cdots\,,
$$
and $0{\scriptsize \bullet} d_1d_2d_3\cdots $ is a different $(-\beta)$-representation of $\ell_{\beta}$.
Nevertheless, the sequence $d_{-\beta}(\ell_\beta)=d_1d_2d_3\cdots$ is of great importance in describing the admissible digit strings.
Before stating Theorem~\ref{thm:admissibility} taken from~\cite{ItoSadahiro}, recall the alternate order on sequences of digits.
We define
$$
v_1v_2v_3\cdots \prec_{\hbox{\tiny alt}} w_1w_2w_3\cdots \quad\hbox{ if }\ \left\{\begin{array}{ll}
v_j<w_j &\hbox{when $j$ is even,}\\
v_j>w_j &\hbox{when $j$ is odd,}
\end{array}\right.\ \hbox{where}\ j:=\min\{i\mid v_i\neq w_i\}\,.
$$
Recall that such ordering is well known from continued fractions. The following theorem is a slight modification of the result of~\cite{ItoSadahiro}, taking into account the requirement~\eqref{eq:ekvivalence}.

\begin{thm}\label{thm:admissibility}
Let $\beta>1$. Let $d_{-\beta}(\ell_\beta)=d_1d_2d_3\cdots$ Define
$$
d_{-\beta}^*(r_\beta)=\left\{
 \begin{array}{ll}
 \big(0d_1\cdots d_{m-1}(d_{m}-1)\big)^\omega &\hbox{if } d_{-\beta}(\ell_\beta)=(d_1\cdots d_m)^\omega,\  \hbox{and $m$ is odd}\,,\\
 0d_1d_2d_3\cdots &\hbox{otherwise.}
 \end{array}\right.
$$
Then the digit sequence $x_kx_{k-1}x_{k-2}\cdots$ is admissible as the $(-\beta)$-expansion of a real number $x$ if and only if each of the suffixes $u$ of the sequence $0x_kx_{k-1}x_{k-2}$ satisfies
$$
d_{-\beta}(\ell_\beta) \preceq_{\hbox{\tiny alt}} u \prec_{\hbox{\tiny alt}} d_{-\beta}^*(r_\beta)\,.
$$
\end{thm}

%

In this paper we focus on the set ${\rm Fin}(-\beta)$ of real numbers with finite number of non-zero digits in their $(-\beta)$-expansion. An important subset of ${\rm Fin}(-\beta)$ is given by the set of so-called $(-\beta)$-integers, denoted by $\Z_{-\beta}$,
$$
\Z_{-\beta}=\{x\in\R\mid \langle x\rangle_{-\beta}=x_kx_{k-1}\cdots x_1x_0{\scriptsize \bullet}\}\,.
$$
In this notation,
$$
{\rm Fin}(-\beta) = \bigcup_{n\in\N}\frac1{(-\beta)^n}\;\Z_{-\beta}\,.
$$
In analogy with the R\'enyi expansions, we study the maximal length of fractional part arising in arithmetical operations, i.e. we investigate
\begin{align*}
L_\oplus(-\beta)&=\min\{l\in\N\mid \forall\,x,y\in\Z_{-\beta},\ x+y\in{\rm Fin}(-\beta)\Rightarrow x+y\in(-\beta)^{-l}\Z_{-\beta}\}\,,\\
L_\otimes(-\beta)&=\min\{l\in\N\mid \forall\,x,y\in\Z_{-\beta},\ x\cdot y\in{\rm Fin}(-\beta)\Rightarrow x\cdot y\in(-\beta)^{-l}\Z_{-\beta}\}\,.
\end{align*}

\section{Arithmetics on $(-\beta)$-expansions}\label{sec:-beta}

The biggest difference between properties of R\'enyi $\beta$-expansions and Ito-Sadahiro $(-\beta)$-expansions
is observed on the set of real numbers which have in their expansion only finitely many non-zero digits. While ${\rm Fin}(\beta)$
is dense in $\R$ for every base $\beta>1$, the set ${\rm Fin}(-\beta)$ may sometimes contain only the point $0$. (Note that $T_{-\beta}(0)=0$ and therefore the string $0^\omega$ is always admissible.)

\begin{thm}\label{thm:fintriv}
Let $\beta>1$. Then ${\rm Fin}(-\beta)=\{0\}$ if and only if $\beta<\frac12(1+\sqrt5)$.
\end{thm}

\pfz
If $x=-\frac1\beta\in\big[-\frac{\beta}{\beta+1},\frac{1}{\beta+1}\big)$, then $\lfloor-\beta x+\frac{\beta}{\beta+1}\rfloor=1$
and $T_{-\beta}(x)=0$. Therefore $d_{-\beta}(x)=10^\omega$, and thus $-\frac1\beta\in{\rm Fin}(-\beta)$.
On the other hand, if $a_{k-1}\cdots a_1a_00^\omega$ is admissible, and at least one digit is non-zero, then Theorem~\ref{thm:admissibility}
implies that also $10^\omega$ is admissible, and thus $-\frac1\beta\in\big[-\frac{\beta}{\beta+1},\frac{1}{\beta+1}\big)$.

By that, we have shown that ${\rm Fin}(-\beta)\neq\{0\}$ if and only if $x=-\frac1\beta\in\big[-\frac{\beta}{\beta+1},\frac{1}{\beta+1}\big)$, which, in turn, is equivalent to the fact that
$$
-\frac{\beta}{\beta+1} \leq -\frac1{\beta}\,.
$$
This is satisfied if and only if $\beta\geq \frac12(1+\sqrt5)$.
\pfk

The proof of the previous theorem implies that if $\beta\geq \tau$, then $1\in{\rm Fin}(-\beta)$. We have $\langle 1\rangle_{-\beta}=1\bullet$ for $\beta>\tau$
and $\langle 1\rangle_{-\beta}=110\bullet$ for $\beta=\tau$ (note that $1=-\tau\ell_{\tau}$).
Our aim is to study for which $\beta\geq\frac12(1+\sqrt5)$ the set ${\rm Fin}(-\beta)$ is a ring. A necessary condition is that $-1\in{\rm Fin}(-\beta)$,
i.e.\ $-1=\sum_{i=k}^{n}a_i(-\beta)^i$, $k,n\in\Z$, $k\leq n$. From this we can derive the following statement.

\begin{lem}\label{l:Y}
Let $\beta>1$ be such that $-1\in{\rm Fin}(-\beta)$. Then $\beta$
is an algebraic number without negative conjugates.
\end{lem}

Assuming that not only $-1$, but all negative integers have finite $(-\beta)$-expansions, we obtain a stronger necessary condition for ${\rm Fin}(-\beta)$
to be a ring. In order to show this, we recall a lemma from~\cite{ADMP}.

\begin{lem}\label{l:intervaly}
Let $a_1a_{2}a_{3}\cdots$ be a
$(-\beta)$-admissible digit string with $a_1\neq 0$. For fixed $k\in\Z$, denote
$$
z=\sum_{i=1}^\infty a_i(-\beta)^{k-i}\,.
$$
Then
$$
z\in\Big[\frac{\beta^{k-1}}{\beta+1},\frac{\beta^{k+1}}{\beta+1}\Big]\,,
\ \hbox{for $k$ odd}\,,\qquad\hbox{ and }\qquad
z\in\Big[-\frac{\beta^{k+1}}{\beta+1},-\frac{\beta^{k-1}}{\beta+1}\Big]\,,
\ \hbox{for $k$ even.}
$$
\end{lem}

\begin{prop}\label{p:finokruh}
Let $\Z^-\subset{\rm Fin}(-\beta)$. Then $\beta$ is a Pisot or a Salem number.
\end{prop}

\pfz
If $\beta$ is an integer, the statement is obvious. Let $\beta\notin\N$.
Theorem~\ref{thm:fintriv} implies that $\beta\geq \tau$. If $\beta=\tau$, the proof is finished.
Consider $\beta>\tau$, i.e. $\beta^2>\beta+1$. First we show that $\beta$ is an algebraic integer.
Consider the negative integer $-\lfloor\beta^{2k+1}\rfloor$ for a fixed $k\in\N$.
Using $\beta^2>\beta+1$ we easily verify that
$$
\frac{\beta^{2k+1}}{\beta+1} < \beta^{2k+1}-1<\beta^{2k+1}<\frac{\beta^{2k+3}}{\beta+1}\,,
$$
and thus $-\lfloor\beta^{(2k+1)}\rfloor\in\Big(-\frac{\beta^{2k+3}}{\beta+1},-\frac{\beta^{2k+1}}{\beta+1}\Big)$.
As $-\lfloor\beta^{2k+1}\rfloor = (-\beta)^{2k+1}+ \varepsilon_\beta$, with $\varepsilon_\beta\in[0,1)$, we can use
Lemma~\ref{l:intervaly} to obtain the $(-\beta)$-expansion of $-\lfloor\beta^{2k+1}\rfloor$  in the form
\begin{equation}\label{eq:2}
-\lfloor\beta^{2k+1}\rfloor = (-\beta)^{2k+1}+ \underbrace{a_0+\frac{a_1}{(-\beta)}+\cdots +\frac{a_l}{(-\beta)^l}}_{\varepsilon_\beta}\,.
\end{equation}
This means that $\beta$ is a root of a monic polynomial with integer coefficients, i.e. $\beta$ is an algebraic integer. In order to show by contradiction
that $\beta$ is a Pisot or Salem number, suppose that $\beta$ has a conjugate $\gamma$ such that $|\gamma|>1$. For $\gamma$, we have
\begin{equation}\label{eq:3}
-\lfloor\beta^{2k+1}\rfloor = (-\gamma)^{2k+1}+
a_0+\frac{a_1}{(-\gamma)}+\cdots +\frac{a_l}{(-\gamma)^l}\,.
\end{equation}
Set $M=\lfloor\beta\rfloor\frac{|\gamma|}{|\gamma-1|}$. Certainly,
there exists a $k\in\N$ such that
\begin{equation}\label{eq:4}
\big|(-\beta)^{2k+1}-(-\gamma)^{2k+1}\big|>1+M\,.
\end{equation}
Since $a_i\in\{0,1,\dots,\lfloor\beta\rfloor\}$, we have
$$
|\varepsilon_\beta|\leq1\quad\hbox{and}\quad|\varepsilon_\gamma|\leq
M\,,\quad\hbox{where}\quad a_0+\frac{a_1}{(-\gamma)}+\cdots
+\frac{a_l}{(-\gamma)^l}\,.
$$
Subtracting~\eqref{eq:2} from~\eqref{eq:3} and comparing
with~\eqref{eq:4} we obtain
$$
1+M<\big|(-\beta)^{2k+1}-(-\gamma)^{2k+1}\big| =
|\varepsilon_\beta-\varepsilon_\gamma|\leq
|\varepsilon_\beta|+|\varepsilon_\gamma|\leq 1+M\,,
$$
which is a contradiction. \pfk

In the following section we study $(-\beta)$-expansions for
$\beta$ a quadratic Pisot number and show that there exist bases
$\beta$ for which ${\rm Fin}(-\beta)$ is a ring. Simpler question
is, whether the set $\Z_{-\beta}$ may be closed under addition for
some $\beta$. It is not difficult to show that if $\beta\in\N$, then $\Z_{-\beta}=\Z$ and
the answer about arithmetical operations is obvious.

\begin{prop} Let $\beta\geq\frac12(1+\sqrt5)$. Then $\Z_{-\beta}$
is a ring if and only if $\beta\in\N$.
\end{prop}

\pfz Let $\beta\notin\N$ and assume that $\Z_{-\beta}$ is a ring.
We have $\lfloor\beta^2\rfloor\in\Z_{-\beta}$. Since
$$
\frac{\beta^2}{\beta+1}<\lfloor\beta^2\rfloor<\frac{\beta^4}{\beta+1}\,,
$$
according to Lemma~\ref{l:intervaly}, the $(-\beta)$-expansion of
$\lfloor\beta^2\rfloor$ is of the form
\begin{equation}\label{eq:5}
\lfloor\beta^2\rfloor =
x_2(-\beta)^2+x_1(-\beta)+x_0\,,\qquad\hbox{where }\
x_0,x_1,x_2\in\{0,1,\dots,\lfloor\beta\rfloor\}\,,\ x_2\geq 1\,.
\end{equation}
This, in turn, means that $\beta$ is a quadratic number.

The assumption $\Z_{-\beta}$ is a ring implies that also ${\rm
Fin}(-\beta)$ is a ring and by Proposition~\ref{p:finokruh},
$\beta$ is either a Pisot or a Salem number. Since there are no
quadratic Salem numbers, $\beta$ is a quadratic Pisot number,
i.e.\ $\beta$ satisfies
\begin{align*}
x^2&=mx+n\,,\qquad\hbox{where }\ m,n\in\N,\ m\geq n,\ \hbox{ or}\\
x^2&=mx-n\,,\qquad\hbox{where }\ m,n\in\N,\ m\geq n+2\geq 3.
\end{align*}
Since the conjugate $\beta'$ of $\beta$ satisfying the equation
$x^2=mx+n$ is negative, by Lemma~\ref{l:Y}, it suffices to
consider the case $\beta^2=m\beta-n$. As
$\lfloor\beta\rfloor=m-1$, the digits $x_0,x_1,x_2$ belong to
$\{0,1,\dots,m-1\}$. From~\eqref{eq:5}, we obtain
$$
\lfloor\beta^2\rfloor = x_2(m\beta-n)-x_1\beta+x_0\,.
$$
Since numbers $\beta,1$ are linearly independent over $\Q$, we
must have $x_2m-x_1=0$, which together with the condition $x_2\geq
1$ and $x_1\leq m-1$ gives a contradiction. \pfk

Determination of values $L_\oplus(-\beta)$, $L_\otimes(-\beta)$ is
in general complicated. For some algebraic numbers $\beta$ with at
least one conjugate $\beta'$ in modulus smaller than 1, we have an
easy estimate which uses the fact that the set $\{|z'|\mid
\Z_{-\beta}\}$ is bounded. Here the notation $z'$ for a number
$z\in\Q(\beta)$ stands for the image of $z$ under the field
isomorphism $\Q(\beta)\mapsto\Q(\beta')$.

\begin{thm}\label{thm:HK}
Let $\beta$ be an algebraic number, and let $\beta'$ be one of
its conjugates satisfying $|\beta'|<1$. Denote
\begin{equation}\label{eq:definiceHK}
\begin{aligned}
H&:=\sup\{|z'|\mid z\in\Z_{-\beta}\}\,,\\
K&:=\inf\{|z'|\mid
z\in\Z_{-\beta}\setminus(-\beta)\Z_{-\beta}\}\,.
\end{aligned}
\end{equation}
If $K>0$, then
\begin{equation}\label{eq:odhadyHK}
\frac1{|\beta'|^{L_\oplus}}\leq\frac{2H}{K} \quad\hbox{and}\quad
\frac1{|\beta'|^{L_\otimes}}\leq \frac{H^2}{K}\,.
\end{equation}
Moreover, if the supremum or infimum in~\eqref{eq:definiceHK} is not reached, then strict inequality holds in both of~\eqref{eq:odhadyHK}.
\end{thm}

\pfz Let $\langle x\pm y\rangle_{-\beta}=z_k\cdots z_1z_0\bullet
z_{-1}\cdots z_{-l}$, for some $x,y\in\Z_{-\beta}$, where
$z_{-l}\neq 0$, i.e.
$$
(x\pm y)(-\beta)^l \in\Z_{-\beta}\setminus(-\beta)\Z_{-\beta}\,.
$$
Thus
$$
K\leq \big|(-\beta)^l(x'\pm y')\big|\leq
|\beta'|^l\big(|x'|+|y'|\big) \leq |\beta'|^l\cdot 2H\,,
$$
which implies
$$
\frac1{|\beta'|^l}\leq \frac{2H}{K}\,.
$$
The statement for $L_\oplus$ is now simple to see. Similarly, we
derive the upper bound for $L_\otimes$.
 \pfk


\section{Quadratic bases in $(-\beta)$-expansions}\label{sec:kvadratic}

The aim of this section is mainly to prove the conjecture of Ito
and Sadahiro~\cite{ItoSadahiro}, which says that for a Pisot
number $\beta$, root of $x^2-mx+n$, $m\geq n+2\geq 3$, the set
${\rm Fin}(-\beta)$ is a ring. For the other class of quadratic
Pisot numbers, roots of $x^2-mx-n$, $m\geq n\geq 1$, this is not
valid. For quadratic Pisot units $\beta$ we obtain values of
$L_\oplus(-\beta)$, $L_\otimes(-\beta)$.

Let us first show that the only quadratic numbers for which
$d_{-\beta}(l_\beta)$ is eventually periodic are quadratic Pisot
numbers. In this, we have an analogue to the case of R\'enyi
$\beta$-expansions, as proved by Bassino~\cite{Bassino}.

\begin{prop}
Let $\beta>1$ be a quadratic number with eventually periodic
$d_{-\beta}(l_\beta)$. Then $\beta$ is a quadratic Pisot number.
\end{prop}

\pfz Let first $d_{-\beta}(l_\beta)=d_1d_2\cdots d_k$ with
$d_k\neq 0$. Then
$$
-\frac{\beta}{\beta+1} =
\frac{d_1}{(-\beta)}+\frac{d_2}{(-\beta)^2}+\cdots+
\frac{d_k}{(-\beta)^k}\,,
$$
which implies
$$
(-\beta)^{k+1} = d_1(-\beta)^{k-1}(\beta+1) +
d_2(-\beta)^{k-2}(\beta+1) + \cdots + d_k(\beta+1)\,.
$$
Therefore $\beta$ is a root of a monic polynomial, say $P(x)$, of
degree $k+1$ with integer coefficients. Thus $\beta$ is a
quadratic integer. Its minimal polynomial is of the form
$Q(x)=x^2-mx-n$, $m,n\in\Z$. As $Q$ divides $P$, necessarily $n$
divides the constant term $P(0)$ of $P$, which satisfies
$|P(0)|=d_k$. Since $d_k$ is a digit of a $(-\beta)$-expansion, we
have $d_k\leq \lfloor\beta\rfloor$, and hence $|n|\leq
\lfloor\beta\rfloor$.

The absolute coefficient of a quadratic polynomial is the product
of $\beta$ and its conjugate $\beta'$, i.e. $|\beta\beta'|=|n|\leq
\lfloor\beta\rfloor<\beta$. This implies that $|\beta'|<1$ and
thus $\beta$ is a Pisot number.

Consider now the situation $d_{-\beta}(l_\beta)=d_1d_2\cdots
d_k(d_{k+1}\cdots d_{k+p})^\omega$. If $k$, $p$ are as small as
possible, we have $d_k\neq d_{k+p}$. Equality
$$
-\frac{\beta}{\beta+1} =
\frac{d_1}{(-\beta)}+\frac{d_2}{(-\beta)^2}+\cdots+
\frac{d_k}{(-\beta)^k} +
\Big(\frac{d_{k+1}}{(-\beta)^{k+1}}+\cdots+
\frac{d_{k+p}}{(-\beta)^{k+p}}\Big)\Big(1+\frac1{(-\beta)^p}+\frac1{(-\beta)^{2p}}+\cdots\Big)\,,
$$
implies that $\beta$ is a root of a monic polynomial of degree
$m+p+1$ with integer coefficients. Moreover, the absolute
coefficient of this polynomial is in modulus equal to
$|d_{k+p}-d_k|$. Since $|d_{k+p}-d_k|\leq \lfloor\beta\rfloor$, we
can derive as before that $|\beta'|<1$.
 \pfk

As a consequence of the above proposition, quadratic numbers $\beta$ with eventually
periodic expansion $d_{-\beta}(l_\beta)$ are precisely the quadratic Pisot numbers, where we have
\begin{align}\label{eq:rozvojekvadrat}
d_{-\beta}(l_\beta) &= \big((m-1)n\big)^\omega &\hbox{for }
\beta^2&=m\beta-n,\quad m-2\geq n\geq 1\,,\\[1mm]\label{eq:rozvojekvadrat2}
d_{-\beta}(l_\beta) &=
m(m-n)^\omega &\hbox{for } \beta^2&=m\beta+n,\quad m\geq n\geq 1\,,
\end{align}
see~\cite{ItoSadahiro}. In the following we study arithmetics for the two classes of quadratic Pisot numbers separately.

\subsection{Case $\boldsymbol{\beta^2=m\beta-n}$, $\boldsymbol{m-2\geq n\geq 1}$}

Let us now focus on the set ${\rm Fin}(-\beta)$ for quadratic Pisot numbers $\beta$,
solutions to $x^2=mx-n$, $m,n\in\N$, $m-2\geq n\geq 1$. For such $\beta$, we have $d_{-\beta}(l_\beta)=\big((m-1)n\big)^\omega$.
and the conjugate of $\beta$ satisfies $\beta'\in(0,1)$. Since $\beta>\tau$, the set  ${\rm Fin}(-\beta)$ is non-trivial, it contains
1. A necessary condition so that ${\rm Fin}(-\beta)$ is closed under addition is that the following implication holds for all $x$,
$$
x\in{\rm Fin}(-\beta)\quad\Rightarrow\quad x+1\in{\rm Fin}(-\beta)\,.
$$
Since we have $(-\beta){\rm Fin}(-\beta)={\rm Fin}(-\beta)$, this condition is also sufficient. Of course, closedness under
addition implies closedness under multiplication. Therefore, if we verify that $-1\in{\rm Fin}(-\beta)$, then closedness under addition
already implies that ${\rm Fin}(-\beta)$ is a ring. For the considered bases $\beta$ we have
$$
\langle -1\rangle_{-\beta}=1(m-1)\bullet n\,.
$$
As a consequence, in order to show that ${\rm Fin}(-\beta)$ is a ring, it suffices to verify the validity of the following statement.

\begin{lem}\label{l:kvadratm-n}
Let $\beta>1$ satisfy $\beta^2=m\beta-n$, for $m,n\in\N$, $m-2\geq n\geq 1$. Then
$x+1\in{\rm Fin}(-\beta)$ for every $x\in{\rm Fin}(-\beta)$.
\end{lem}

Before proving the lemma, let us first describe the digit strings that are admissible as $(-\beta)$-expansions.
The following is a simple consequence of Theorem~\eqref{thm:admissibility}, taking into account expression~\eqref{eq:rozvojekvadrat}.

\begin{lem}\label{l:m-nzakazane}
Let $\beta$ be the larger root of $x^2-mx+n$, $m-2\geq n\geq 1$. A digit string $x_kx_{k-1}x_{k-2}\cdots$ with finitely many non-zero digits
is admissible as $(-\beta)$-expansion if and only if $x_i\in\{0,1,\dots,m-1\}$ and $x_{i}=m-1 \ \Rightarrow \ x_{i-1}\geq n$.
\end{lem}

\pfz[Proof of Lemma~\ref{l:kvadratm-n}]
Let us first realize that since $\beta$ is the root of $x^2-mx+n$, we have the following representations of 0,
$$
1\ m\ n \ \bullet= \overline{1}\ \overline{m}\ \overline{n} \ \bullet= 0\,,
$$
where $\overline{A}$ is a compact form of $-A$.
By repeated application of this relation, we obtain for every $k\in\N$,
\begin{equation}\label{eq:zero}
\begin{array}{ccccccccccc}
   & 1 & (m-1) & [\ol{(m-n-1)} & (m-n-1)]^k & \ol{(m-n)} & \ol{n} &\bullet & = 0 \\
  1 & (m-1) & \ol{(m-n-1)} & [(m-n-1) & \ol{(m-n-1)}]^k & (m-n) & n &\bullet & =0\\
\end{array}
\end{equation}

Adding 1 to a number $x$ written as an admissible digit string may result in a non-admissible digit string, which, nevertheless, represents the number $x+1$. We show that $x+1$ belongs to ${\rm Fin}(-\beta)$ by providing its finite $(-\beta)$-expansion. In order to see that the two strings represent the same number, one can verify that the second one is obtained from the first one by adding digit-wise a zero which is in the form~\eqref{eq:zero}.
We give a list of cases. One verifies by inspection that the list contains all cases of non-admissible strings that arise from admissible ones by adding 1.

According to Lemma~\ref{l:m-nzakazane}, a digit string may be non-admissible by breaking one of the two conditions, namely, either it is not over the alphabet
$\{0,1,\dots,m-1\}$, or it contains the subsequence $(m-1)A$, where $A\leq n-1$.

{\renewcommand{\arraystretch}{1.5}
\noindent{\bf Case 1.} Consider an $x\in{\rm Fin}(\beta)$ such that its
$(-\beta)$-expansion has digit $m-1$ at $(-\beta)^0$. Then
necessarily the digit at position $(-\beta)^{-1}$, denote it by
$C$, is at least $n$. Find $k\in\{0,1,2,\dots\}$ such that we have
a representation of $x+1$ in the form
\begin{equation}\label{e:4}
x+1 \ = \ \cdots\ A\ B\ [(m\!-\!1)\ n]^k\ m\ {\scriptsize
\bullet}\ C\ \cdots\qquad \hbox{ where the string } \ A\ B \neq (m\!-\!1)\
n\,.
\end{equation}

\noindent{\bf Case 1.1.} First take $B=0$. To the representation~\eqref{e:4} of
the number $x+1$ we add digit-wise a representation of $0$,
$$
 \begin{array}{cccccccc@{\quad {\scriptsize \bullet}\quad }cc}
  x+1&=&\cdots & A & 0 & [(m\!-\!1) & n]^k & m & C &\cdots \\
  0&=& &   1 & (m\!-\!1) & [\ol{(m\!-\!n\!-\!1)} & (m\!-\!n\!-\!1)]^k & \ol{(m\!-\!n)} & \ol{n}& \\
  \hline
  x+1&=&\cdots & (A\!+\!1)& (m\!-\!1)& [n & (m\!-\!1)]^k & n & C\!-\!n&\cdots
  \end{array}
$$
Since $B=0$ in the $(-\beta)$-expansion of $x$, we necessarily
have $A\leq m-2$. Therefore also the resulting representation of
$x+1$ is admissible as $(-\beta)$-expansion of $x+1$.

The case $B\geq1$ is divided into two.

\noindent{\bf Case 1.2.} Let $B\geq 1$ and $k=0$. Again, we add to the
non-admissible representation of $x+1$ in the form~\eqref{e:4} a
suitable representation of $0$,
$$
 \begin{array}{cccccc@{\quad {\scriptsize \bullet}\quad }cc}
  x+1&=&\cdots & A & B & m & C &\cdots \\
  0&=& &   & \ol{1} & \ol{m} & \ol{n} \\
  \hline
  x+1&=&\cdots & A& (B\!-\!1)& 0& C\!-\!n&\cdots
  \end{array}
$$
Since $A\ B\neq (m-1)\ n$, the resulting representation of $x+1$
is the $(-\beta)$-expansion of $x+1$.

\noindent{\bf Case 1.3.} Let $B\geq 1$ and $k\geq 1$. In this case we rewrite
{\small
$$
 \begin{array}{cccccccccc@{\quad {\scriptsize \bullet}\quad }cc}
  x+1&=&\cdots & A & B&(m\!-\!1)&n & [(m\!-\!1) & n]^{k\!-\!1} & m & C &\cdots \\
  0&=& &   & \ol{1} & \ol{(m\!-\!1)}&(m\!-\!n\!-\!1)&[\ol{(m\!-\!n\!-\!1)} & (m\!-\!n\!-\!1)]^{k\!-\!1} & \ol{(m\!-\!n)} & \ol{n}& \\
  \hline
  x+1&=&\cdots &A & (B\!-\!1)& 0& (m\!-\!1)&[n & (m\!-\!1)]^{k\!-\!1} & n & C\!-\!n&\cdots
  \end{array}
$$}

\noindent{\bf Case 2.} Consider an $x\in{\rm Fin}(\beta)$ such that its
$(-\beta)$-expansion has digit $m-2$ at the position $(-\beta)^0$.
In order that after adding 1 one obtains a non-admissible string,
necessarily the digit $C$ at position $(-\beta)^{-1}$ satisfies
$C\leq n-1$. Denote by $M$ the set of pairs of digits
$$
M:=\big\{\,X\ Y \ \big|\ X\in\{m-n-1,\dots,m-1\},\
Y\in\{0,1,\dots,n-1\}\,\big\}\,.
$$
Then we can find $k,l\in\{0,1,2,\dots\}$ such that
\begin{equation}\label{e:5}
x+1 \ = \ \cdots \ A\ B\ [(m\!-\!1)\ n]^k\ (m\!-\!1)\ {\scriptsize
\bullet}\ C\ X_1\ Y_1 \cdots X_l\ Y_l\ D\ E\ \cdots
\end{equation}
where the string $A\ B \neq (m\!-\!1)\ n$ and the string $D\ E$ does not belong to $M$. Denote by
$p_1$, $p_2$ the $(-\beta)$-integers
$$
\begin{array}{cccccccc}
p_1 &=&  &\ol{1} & \ol{(m\!-\!1)}  &[{(m\!-\!n\!-\!1)} & \ol{(m\!-\!n\!-\!1)}]^k & {\scriptsize \bullet}\\
p_2 &=& 1 & (m\!-\!1) & \ol{(m\!-\!n\!-\!1)} &[{(m\!-\!n\!-\!1)} &
\ol{(m\!-\!n\!-\!1)}]^k &{\scriptsize \bullet}
\end{array}
$$
and by $z_1$, $z_2$ the following numbers with only
$(-\beta)$-fractional part,
$$
\begin{array}{cccccccc}
z_1 &=&  {\scriptsize \bullet} & [{(m\!-\!n\!-\!1)} & \ol{(m\!-\!n\!-\!1)}]^l & (m\!-\!n) & n\\
z_2 &=&  {\scriptsize \bullet} & [{(m\!-\!n\!-\!1)} &
\ol{(m\!-\!n\!-\!1)}]^l & (m\!-\!n\!-\!1)&\ol{(m\!-\!n)} & \ol{n}
\end{array}
$$
Using~\eqref{eq:zero} one can easily see that $p_i+z_j = 0$ for
$i,j\in\{1,2\}$. We shall work separately with the
$(-\beta)$-integer and $(-\beta)$-fractional part of $x+1$.

First consider the $(-\beta)$-fractional part of $x+1$. Recall
that $C\leq n-1$ and the string $D\ E\notin M$.
If $D\leq m-n-2$ or $D=m-n-1$ (the latter implies $E\geq n$), then
we add $z_1$ to the $(-\beta)$-fractional part of $x+1$. We obtain
{\small
$$
\begin{array}{@{\quad{\scriptsize \bullet}\quad}ccccccccc}
 C&X_1&Y_1 &\cdots & X_l & Y_l & D & E & \cdots\\
 {(m\!-\!n\!-\!1)}& \ol{(m\!-\!n\!-\!1)}& {(m\!-\!n\!-\!1)}
  &\cdots & \ol{(m\!-\!n\!-\!1)}& {(m\!-\!n)} & n\\ \hline
 C\!+\!m\!-\!n\!-\!1 &X_1\!-\!(m\!-\!n\!-\!1) & Y_1\!+\!m\!-\!n\!-\!1
 &\cdots &X_l\!-\!(m\!-\!n\!-\!1) & Y_l\!+\!m\!-\!n & D\!+\!n & E
 & \cdots
\end{array}
$$}
The resulting fractional part is an admissible digit string.
If $D\geq m-n$ (which implies $E\geq n$), then we add $z_2$ to the
$(-\beta)$-fractional part of $x+1$. We obtain
{\small
$$
\begin{array}{@{\quad{\scriptsize \bullet}\quad}cccccccccc}
 C&X_1&Y_1 &\cdots & X_l & Y_l & D & E & \cdots\\
 {(m\!-\!n\!-\!1)}& \ol{(m\!-\!n\!-\!1)}& {(m\!-\!n\!-\!1)}
  &\cdots & \ol{(m\!-\!n\!-\!1)}& {(m\!-\!n\!-\!1)}&\ol{(m\!-\!n\!)} & \ol{n}\\ \hline
 C\!+\!m\!-\!n\!-\!1 &X_1\!-\!(m\!-\!n\!-\!1) & Y_1\!+\!m\!-\!n\!-\!1
 &\cdots &X_l\!-\!(m\!-\!n\!-\!1) & Y_l\!+\!m\!-\!n\!-\!1 & D\!-\!m\!+\!n &
 E\!-\!n
 & \cdots
\end{array}
$$}
Again, the resulting string is admissible.

Let us now take the $(-\beta)$-integer part of $x+1$. Recall that
$A\ B\neq (m\!-\!1)\ n$. If $B=0$, then to the $(-\beta)$-integer
part of $x+1$ we add $p_2$, if $B\geq 1$, we add $p_1$. We obtain
$$
\begin{array}{cccccc@{\quad{\scriptsize \bullet}\quad}}
 \cdots& A& 0 & [(m\!-\!1) & n]^k & (m\!-\!1)\\
 &1&(m\!-\!1)& [\ol{(m\!-\!n\!-\!1)} & {(m\!-\!n\!-\!1)}]^k &
 \ol{(m\!-\!n\!-\!1)}\\ \hline
 \cdots& (A\!+\!1)& (m\!-\!1) & [n&(m\!-\!1) ]^k & n
\end{array}
$$
and
$$
\begin{array}{cccccc@{\quad{\scriptsize \bullet}\quad}}
 \cdots& A& B & [(m\!-\!1) & n]^k & (m\!-\!1)\\
 &&\ol{1}&\ol{(m\!-\!1)}& [{(m\!-\!n\!-\!1)} & \ol{(m\!-\!n\!-\!1)}]^k \\ \hline
 \cdots& A & (B\!-\!1)& 0  & [(m\!-\!1)&n ]^k
\end{array}
$$
In both cases, the result is an admissible string with last digit
equal to $n$. Concatenating such a string with an admissible digit
string resulting from the $(-\beta)$-fractional part, we obtain an
admissible digit string. Therefore, we have provided a
prescription to rewrite the original non-admissible representation
of $x+1$ of the form~\eqref{e:5} by adding $0$ in the form
$p_i+z_j$, into the $(-\beta)$-expansion of $x+1$.
This completes the proof.} \pfk

\begin{coro}
Let $\beta>1$ satisfy $\beta^2=m\beta-n$, for $m,n\in\N$, $m-2\geq n\geq 1$. Then ${\rm Fin}(-\beta)$ is a ring.
\end{coro}

Addition of two $(-\beta)$-integers thus always yields a number whose $(-\beta)$-expansion has a finite number of fractional digits.
We can give the upper bound to the length of such a fractional part in case that $\beta$ is a unit.

\begin{ex}\label{ex:quadrunitmistaA}
Let $\beta>1$ satisfy $\beta^2=m\beta-1$, for $m\in\N$, $m\geq 3$. The digits of $(-\beta)$-expansions thus take values in the set
$\{0,1,\dots,m-1\}$, and $(m-1)0$ is a forbidden string. Put $x=m-2$, $y=1$. Obviously $x,y\in\Z_{-\beta}$. In order to find the
$(-\beta)$-expansion of $z=x+y=m-1$, we rewrite
$$
m-1 = (-\beta)^2 + (m-1) (-\beta) + 1 + \frac{m-1}{(-\beta)} + \frac1{(-\beta)^2}\,,
$$
and hence
$$
\langle z\rangle_{-\beta}=1(m-1)1\bullet (m-1)1\in\frac1{(-\beta)^{2}}\Z_{-\beta}\setminus\frac1{(-\beta)}\Z_{-\beta}\,.
$$
From this, we can conclude that $L_\oplus(-\beta)\geq 2$.

Now put $x=(m-2)(1+\beta^2)$, $y=1+\beta^2$. The corresponding $(-\beta)$-expansions are $\langle x\rangle_{-\beta}=(m-2)0(m-2)\bullet$
and $\langle y\rangle_{-\beta}=101\bullet$, i.e.\ $x,y\in\Z_{-\beta}$. One easily verifies that
$$
\begin{aligned}
z=xy &=(m-2)\beta^4
+2(m-2)\beta^2+(m-1)= \\[1mm]
&= (m-1)(-\beta)^4 + (m-1)(-\beta)^3 + (m-2)(-\beta)^2 + (m-2)(-\beta)  + \frac{m-1}{(-\beta)} + \frac1{(-\beta)^2}\,,
\end{aligned}
$$
and hence
$$
\langle z\rangle_{-\beta}=(m-1)(m-1)(m-2)(m-2)0\bullet (m-1)1\in\frac1{(-\beta)^{2}}\Z_{-\beta}\setminus\frac1{(-\beta)}\Z_{-\beta}\,,
$$
and we can conclude that $L_\otimes(-\beta)\geq 2$.
\end{ex}

\begin{thm}
Let $\beta>1$ satisfy $\beta^2=m\beta-1$, for $m\in\N$, $m\geq 3$. Then
$$
L_\oplus(-\beta) = 2 = L_{\otimes}(-\beta)\,.
$$
\end{thm}

\pfz
We use Theorem~\ref{thm:HK}. We first show that
$$
H= \frac{1-\beta'}{{\beta'}(1+\beta')} \quad\hbox{ and }\quad K=\frac{{\beta'}(1-\beta')}{1+\beta'}\,.
$$
Consider $z\in\Z_{-\beta}$ with the expansion $\langle z\rangle_{-\beta} = z_kz_{k-1}\cdots z_1z_0\bullet$, where
$z_i=m-1$ for at least one index $i$. Since the string $(m-1)0$ is forbidden, we have $i\geq 1$ and, moreover, $z_{i-1}\geq 1$.
Since the string $z_k\cdots z_{i+1}(m-1)z_{i-1}\cdots z_1z_0\bullet$ is admissible, also the strings $z_k\cdots z_{i+1}(m-2)z_{i-1}\cdots z_1z_0\bullet$
and $z_k\cdots z_{i+1}(m-2)(z_{i-1}-1)\cdots z_1z_0\bullet$ are admissible, denote by $x$ and $y$ the corresponding $(-\beta)$-integers. For the field conjugates $x'$, $y'$, $z'$ of $x$, $y$, $z$ we have
$$
x'+(-\beta)^i = z' = y' + (-\beta)^i + (-\beta)^{i-1} = y'+(-\beta)^{i-1}(1-\beta')\,.
$$
Since $\beta'\in(0,1)$, the number $z'$ lies between numbers $x'$ and $y'$. From this, we can derive that
$$
H=\sup\big\{|z'|\mid z\in\Z_{-\beta},\ \langle z\rangle_{-\beta} \hbox{ does not contain the digit }(m-1)\big\}\,.
$$
Note that arbitrary string of digits $\{0,1,\dots,m-2\}$ is admissible. For $z\in\Z_{-\beta}$ we can thus write
$$
z'=\sum_{i=0}^kz_i(-\beta')^i\leq \sum_{0\leq 2i\leq k}(m-2)(-\beta')^{2i}\leq \frac{m-2}{1-(\beta')^2}\,.
$$
Similarly, we have
$$
z'=\sum_{i=0}^kz_i(-\beta')^i\geq \sum_{0\leq 2i+1\leq k}(m-2)(-\beta')^{2i+1}\geq -\beta'\frac{m-2}{1-(\beta')^2}\,.
$$
Equality $(\beta')^2=m\beta'-1$ implies $(m-2)\beta'=(\beta'-1)^2$ and therefore
$$
H=\frac{m-2}{1-(\beta')^2}=\frac1{\beta'}\frac{(1-\beta')^2}{1-(\beta')^2} = \frac1{\beta'}\frac{1-\beta'}{1+\beta'}\,,
$$
as we wanted to show.

In order to determine $K$, we study the field conjugate $z'$ for $z\in \Z_{-\beta}\setminus(-\beta)\Z_{-\beta}$. Again, for $z$ whose $(-\beta)$-expansion
 contains the digit $(m-1)$ we can find $x,y\in\Z_{-\beta}\setminus(-\beta)\Z_{-\beta}$ with digits in $\{0,1,\dots,m-2\}$ such that $x'<z'<y'$. The only exception is the case when the last two digits of $z$ equal to $z_1z_0=(m-1)1$. Here we have
 $$
 z'=\sum_{i=0}^kz_i(-\beta')^i > 1+(m-1)(-\beta') + \sum_{i=0}^\infty (m-2)(-\beta')^{2i+1} = 1-\beta' -\beta' H =
 \frac{{\beta'}(1-\beta')}{1+\beta'} = K\,.
 $$
Let us put the computed values of $H,K$ into Theorem~\ref{thm:HK}. With the use of $\beta'=\beta^{-1}$ we obtain that
\begin{align*}
\beta^{L_\oplus} &\leq \frac{2H}{K} = \frac2{(\beta')^2}=2\beta^2<\beta^3\,,\\
\beta^{L_\otimes} &\leq \frac{H^2}{K} = \frac1{(\beta')^3}\frac{1-\beta'}{1+\beta'}=\beta^3\frac{1-\beta'}{1+\beta'}<\beta^3\,.
\end{align*}
Together with Example~\eqref{ex:quadrunitmistaA}, we have the statement of the theorem.
\pfk

\subsection{Case $\boldsymbol{\beta^2=m\beta+n}$, $\boldsymbol{m\geq n\geq 1}$}

Consider for $\beta$ the larger root of $x^2=mx+n$, $m\geq n\geq 1$. Its conjugate $\beta'$ is negative, and therefore by
Lemma~\ref{l:Y} the set ${\rm Fin}(-\beta)$ cannot be a ring. In particular, we have an infinite expansion for the number $-1$,
$$
\langle -1\rangle_{-\beta} = 1m\bullet(m-n+1)^\omega\,.
$$
Nevertheless, this fact does not prevent ${\rm Fin}(-\beta)$ to be closed under addition, as it is shown in~\cite{Vavra}.

Even if by summing and multiplying $(-\beta)$-integers leads sometimes to infinite $(-\beta)$-expansions,
we can estimate the number of arising fractional digits in case that the result of addition and multiplication belongs to
${\rm Fin}(-\beta)$. Theorem~\ref{thm:HK} leads a simple result for bases $\beta$ which are units.

\begin{thm}
Let $\beta>1$ satisfy $\beta^2=m\beta+1$, for $m\in\N$, $m\geq2$. Then
$$
L_\oplus(-\beta) = 1 = L_{\otimes}(-\beta)\,.
$$
\end{thm}

\pfz
First we show that the values $L_\oplus(-\beta) = 1$,  $L_{\otimes}(-\beta)=1$ can be reached.
The digits of $(-\beta)$-expansions take values in the set
$\{0,1,\dots,m\}$. Using~\eqref{thm:admissibility} and~\eqref{eq:rozvojekvadrat2} we derive that forbidden are
the strings $m(m-1)^{2k+1}m$  and $m(m-1)^{2k}A$, where $k\geq 0$ and $A\leq m-2$.
Put $x=-m\beta + m-1$. Obviously $x\in\Z_{-\beta}$. We rewrite
$$
z=x+x=2m(-\beta) + (m-1) = (-\beta)^3+ (m-1)(-\beta)^2 + (m-2)(-\beta)+ m-1 + \frac{1}{(-\beta)}\,,
$$
and hence
$$
\langle z\rangle_{-\beta}=1(m-1)(m-2)(m-1)\bullet 1\in\frac1{(-\beta)}\Z_{-\beta}\setminus\Z_{-\beta}\,.
$$
From this, we can conclude that $L_\oplus(-\beta)\geq 1$.
In order to find the lower bound to $L_\otimes(-\beta)$, it suffices to realize that $y=2\in\Z_{-\beta}$. For, we have
$\langle y \rangle_{-\beta}=2\bullet$ if $m\geq 3$, and $\langle y \rangle_{-\beta}=121\bullet$ if $m= 2$.
We therefore have $xy = z$ as above and $L_\otimes(-\beta)\geq 1$.

In order to find upper bounds on $L_\oplus(-\beta)$, $L_\otimes(-\beta)$, let us compute the values $H$, $K$ for use in Theorem~\ref{thm:HK}.
Consider $z\in\Z_{-\beta}$ with the expansion $\langle z\rangle_{-\beta} = z_kz_{k-1}\cdots z_1z_0\bullet$.
Since the conjugate $\beta'=-\beta^{-1}$ of $\beta$ belongs to $(-1,0)$, we have
$$
0\leq z'=\sum_{i=0}^k z_i(-\beta')^i = \sum_{i=0}^kz_i|\beta'|^i\leq m-1 + \sum_{i=1}^\infty m|\beta'|^i = \frac{m\beta}{\beta-1}-1 = \beta\,,
$$
where we have taken account of the fact that a $(-\beta)$-integer cannot have the digit $m$ at the position $(-\beta)^0$.
Therefore $H=\sup\{|z'|\mid z\in\Z_{-\beta}\}=\beta$ and the supremum is not reached. Similarly, the fact that $\beta'\in(-1,0)$ implies that $K=1$.
We therefore have
\begin{align*}
\beta^{L_\oplus} &< \frac{2H}{K} = 2\beta <\beta^2\,,\\
\beta^{L_\otimes} &< \frac{H^2}{K} = \beta^2\,,
\end{align*}
where we have used that $\beta>2$, which is valid, as we consider $m\geq 2$. This implies the statement of the theorem.
\pfk

Note that we derive the exact values of $L_\oplus(-\beta)$,
$L_{\otimes}(-\beta)$ only for quadratic Pisot units not equal to
the golden ratio. In fact, Theorem~\ref{thm:HK} leads in the case
$\beta=\frac12(1+\sqrt5)$ only to the values $L_\oplus(-\beta)
\leq 3$ and  $L_{\otimes}(-\beta)\leq 3$, but the actual values
are $L_\oplus(-\beta) = 2 = L_{\otimes}(-\beta)$. The proof of
this fact is found in~\cite{Vavra}.

\section*{Acknowledgements}

We acknowledge financial support by the Czech Science Foundation
grant 201/09/0584 and by the grants MSM6840770039 and LC06002 of
the Ministry of Education, Youth, and Sports of the Czech
Republic. The work was also partially supported by the CTU student grant SGS10/085/OHK4/1T/14.



\begin{thebibliography}{99}

\bibitem{akiyamaCubic}
S.~Akiyama, {\it Cubic Pisot Units with finite beta expansions}, in `Algebraic Number Theory and Diophantine Analysis',
ed. by F.Halter-Koch and R.F. Tichy, de Gruyter (2000), 11--26.

\bibitem{akiyamaTiles}
S.~Akiyama, {\it Self affine tiling and Pisot numeration system}, 'Number Theory
and its Applications', ed. by K. Gy¨ory and S. Kanemitsu, 7–-17 Kluwer 1999.

\bibitem{ADMP} P.~Ambro\v z, D.~Dombek, Z.~Mas\'akov\'a, E.~Pelantov\'a,
{\it Numbers with integer expansion in the numeration system with negative base},
preprint 2009, 13pp. {\tt http://arxiv.org/abs/0912.4597}

\bibitem{AFMP} P.~Ambro\v z, C. Frougny, Z.~Mas\'akov\'a, E.~Pelantov\'a,
{\it Arithmetics on number systems with irrational bases}, Bull. Belgian Math. Soc. Simon Stevin {\bf 10} (2003), 641–-659.

\bibitem{Bassino}
F.~Bassino, {\it $\beta$-expansions for cubic Pisot numbers},  5th Latin American Theoretical INformatics Symposium
(LATIN'02), {\bf 2286} LNCS.  Cancun, Mexico.  April,  2002. pp. 141-–152. Springer-Verlag.

\bibitem{Bernat} J. Bernat, {\it Arithmetics in $\beta$-numeration}, Discr. Math. Theor. Comp. Sci.  {\bf 9} (2007), 85--106.

\bibitem{BernatCubic} J. Bernat, {\it Computation of $L_\oplus$ for several cubic Pisot numbers}, Discr. Math. Theor. Comput. Sci. {\bf 9} (2007), 175--193.

\bibitem{BuFrGaKr} \v{C}.~Burd\'\i k, Ch.~Frougny, J.~P.~Gazeau, R.~Krejcar,
{\it Beta-Integers as Natural Counting Systems for Quasicrystals ,}
J.~Phys. A: Math. Gen. {\bf 31} (1998) 6449--6472.

\bibitem{ChiaraFrougny}
Ch.~Frougny, A.~C.~Lai, {\it On negative bases}, Proceedings of DLT 09, Lectures Notes in Computer Science, 5583 (2009).

\bibitem{FruSo} Ch.~Frougny and B.~Solomyak, {\it Finite
$\beta$-expansions,} Ergodic Theory Dynamical Systems {\bf 12}
(1994), 713--723.

\bibitem{GuMaPe} L.~S.~Guimond, Z.~Mas\'akov\'a, E.~Pelantov\'a, {\it Arithmetics of beta-expansions}, Acta Arith.
{\bf 112} (2004), 23–-40.

\bibitem{Hollander}
M.~Hollander, {\it Linear Numeration systems, Finite Beta Expansions, and Discrete
Spectrum of Substitution Dynamical Systems}, Ph.D. thesis, University
of Washington, 1996.

\bibitem{ItoSadahiro}
  S.~Ito, T.~Sadahiro,
  {\it $(-\beta)$-expansions of real numbers},
  Integers  {\bf 9} (2009), 239--259.


\bibitem{Vavra} Z. Mas\'akov\'a, T. V\'avra, {\it Arithmetics in numeration systems with negative quadratic base},
Kybernetika {\bf 47} (2011), (in press).

\bibitem{Parry} W.~Parry, {\it On the $\beta$-expansions of real numbers,}
Acta Math. Acad. Sci. Hung. {\bf 11} (1960), 401--416.

\bibitem{Renyi}
A. R\'enyi, {\it Representations for real numbers and their
ergodic properties,} Acta Math. Acad. Sci. Hung. {\bf 8} (1957),
477--493.



\end{thebibliography}
\end{document}